\documentclass[12pt]{elsarticle}
\usepackage[utf8]{inputenc}
\usepackage{amsmath}
\usepackage{bbold}
\usepackage{latexsym}
\usepackage{epsfig}
\usepackage{psfrag}
\usepackage{enumerate}
\usepackage[utf8]{inputenc}
\usepackage{newtxtext}
\usepackage{xcolor}
\usepackage{graphicx}
\makeatletter
\def\ps@pprintTitle{%
 \let\@oddhead\@empty
 \let\@evenhead\@empty
 \def\@oddfoot{\centerline{\thepage}}%
 \let\@evenfoot\@oddfoot}
\makeatother

\newtheorem{thm}{Theorem}[section]

\newtheorem{conjecture}{Conjecture}[section]

\usepackage{amssymb}
\begin{document}
\begin{frontmatter}

\title{ A simple proof of  the Gr\"unbaum conjecture.}

 \author[label1]{Beata~Der\c{e}gowska\footnote{B.D. is partially supported by National Science Center (NCN) grant no. 2021/05/X/ST1/01212. For the purpose of Open Access, the author has applied a CC-BY public copyright licence to any Author Accepted Manuscript (AAM) version arising from this submission. }}
 
 \author[label4]{Barbara~Lewandowska}
 \address[label1]{Institute of Mathematics\\
Pedagogical University of Krakow, Podchorazych~2, Krakow, 30-084, Poland}
\address[label4]{Faculty of Mathematics and Computer Science\\
Jagiellonian University, Lojasiewicza~6, Krakow, 30-048, Poland}

%% \author{llll}

%% \affiliation{organization={},%Department and Organization

%\date{June 2022}
\begin{abstract}
Let $\lambda_\mathbb{K}(m)$ denote the maximal absolute projection constant over the subspaces of dimension $m$. Apart from the trivial case for $ m=1$, the only known value of $\lambda_\mathbb{K}(m)$ is for $ m=2$ and $\mathbb{K}=\mathbb{R}.$  In 1960, B.Grünbaum conjectured that $\lambda_\mathbb{R}(2)=\frac{4}{3}$ and, in 2010, B. Chalmers and G. Lewicki proved it. In 2019, G. Basso delivered the alternative proof of this conjecture. Both proofs are quite complicated, and there was a strong belief that providing an exact value for $\lambda_\mathbb{K}(m)$ in other cases would be a difficult task. In our paper, we present an upper bound of the value $\lambda_\mathbb{K}(m)$, which becomes an exact value for numerous cases. This bound was first stated in [H. K\"{o}nig,  N. Tomczak-Jaegermann,
{\it Norms of minimal projections,} 
J. Funct. Anal. 119/2 (1994),  253--280], but with an erroneous proof. The crucial idea of our proof will be an application of some results from the articles [B. Bukh, C. Cox, {\it Nearly orthogonal vectors and small antipodal spherical codes,} Isr. J. Math. 238, 359–388 (2020)] and [G. Basso, {\it Computation of maximal projection constants,} J. Funct. Anal. 277/10 (2019), 3560--3585.], for which simplified proofs will be given.

\end{abstract}
\end{frontmatter}
\section{Introduction}
Let $X$ be a Banach space over $\mathbb{K}$ (where $\mathbb{K}=\mathbb{R}$ or $\mathbb{K}=\mathbb{C}$) and let $Y\subset X$ be a finite-dimensional subspace.
Let $\mathcal{P}(X,Y)$ denote the set of all linear and continuous projections from $X$ onto $Y$,
recalling that an operator $P \colon X \rightarrow Y$ is called a \textit{projection} onto $Y$ if $P|_Y={\rm Id}_Y.$ 
We define the \textit{relative projection constant} of $Y$ by
\begin{equation*}
\lambda_{\mathbb{K}}(Y,X) :=\inf\lbrace\|P\|:\;P\in\mathcal{P}(X,Y)\rbrace
\end{equation*}
and the \textit{absolute projection constant} of $Y$ by
\begin{equation}
\label{DefMAPC}
\lambda_{\mathbb{K}}(Y) :=\sup\lbrace\lambda(Y,X):Y\subset X\rbrace,
\end{equation}
and finally  the  \textit{maximal absolute projection  constant},  by
\begin{equation*}
\lambda_{\mathbb{K}}(m) :=\sup \lbrace\lambda(Y):\; \dim(Y)=m \rbrace.
\end{equation*}
To calculate this value, it suffices to take the supremum over finite-dimensional $l_\infty(\mathbb{K})$ superspaces (see e.g. \cite[III.B.5]{W}). Therefore, the latter can be defined as a supremum over \textit{maximal relative projection constants} for $N \ge m$ 
\begin{equation*}
\lambda_{\mathbb{K}}(m,N):=\sup\lbrace \lambda(Y, l_\infty^{(N)}(\mathbb{K})):\; \dim(Y)=m \textrm{ and } Y\subset l_\infty^{(N)}(\mathbb{K})\rbrace.
\end{equation*}   
Our estimation of $\lambda_{\mathbb{K}}(m,N)$
will rely on the following result proved in \cite[combining Theorem 2.2 and Theorem 2.1]{CLe}  (the simplified version of the proof can be found in \cite[ Appendix A]{FS}).

\begin{thm}\label{lammbda}
For integers $N \ge m$, we have
\begin{equation*}
    \lambda_{\mathbb{K}}(m,N)=\max\bigg\lbrace \sum_{i,j=1}^N t_it_j|U^* U|_{ij}:t\in\mathbb{R}_+^N,\;\|t\|=1,U\in \mathbb{K}^{m\times N},\; UU^*={\rm I}_m \bigg\rbrace .
\end{equation*}
\end{thm}
Since computing this value is difficult, the literature deals with its lower bound called the {\it quasimaximal relative projection constant}, which arises when choosing a vector $t$ with equal coordinates. To be more precise, for $N\geq m$
\begin{equation}\label{mi}
    \mu_\mathbb{K}(m,N):=\max\bigg\lbrace \frac{1}{N}\sum_{i,j=1}^{N}|U^* U|_{ij}:U\in \mathbb{K}^{m\times N},\; UU^*={\rm I}_m \bigg\rbrace .
\end{equation}
Analogously as for the maximal relative projection constant, we define the {\it quasimaximal absolute projection constant} by
\begin{equation}
    \mu_\mathbb{K}(m)=\sup\{\mu_\mathbb{K}(m,N): N\geq m\}.
\end{equation}
In fact, it appears that $\mu_\mathbb{K}(m)$ is an equivalent definition of $\lambda_\mathbb{K}(m).$ In the real case, it was proved by Basso \cite[Proof of Theorem 1.2]{B}. In our paper, we give an elementary proof of this fact, which is also valid in the complex case (Theorem \ref{basso}). The best known bound for the maximal relative projection constant was given in \cite{KLL}.  We present it in the form stated in \cite [ Theorem 5 ]{ FS}, where also an easier proof of this result was provided.

\begin{thm}\label{SF}

For integers $N\geq m$, the maximal relative projection constant $\lambda_{\mathbb{K}}(m,N)$ is upper bounded by
$$
\delta_{m,N} := \frac{m}{N} \left( 1 + \sqrt{\frac{(N-1)(N-m)}{m}} \right).
$$
Moreover, the following properties are equivalent:
\begin{enumerate}[i)]
\item There is an equiangular tight frame consisting of $N$ vectors in $\mathbb{K}^m,$
\item $\mu_\mathbb{K}(m,N)=\tfrac{m}{N}\left(1 +\sqrt{\tfrac{(N-1)(N-m)}{m}} \right),$
\item $\lambda_\mathbb{K}(m,N)=\tfrac{m}{N}\left(1 +\sqrt{\tfrac{(N-1)(N-m)}{m}} \right).$
\end{enumerate}
\end{thm}
Recall that a system of vectors $(u_1,\dots, u_N)$ in $\mathbb{K}^m$ is called a {\it tight frame} if there exists a constant $\alpha >0$ such that one of the following equivalent conditions holds:
\begin{itemize}
    \item  $\|x\|^2=\alpha\sum_{k=1}^{N}|\langle x, u_k \rangle|^2$  \; for all $x\in \mathbb{K}^m.$ 
     \item  $x=\alpha\sum_{k=1}^{N}\langle x, u_k \rangle u_k$  \; for all $x\in \mathbb{K}^m.$ 
    \item $UU^*=\frac{1}{\alpha}{\rm I_m}$, where $U$ is the matrix with columns $u_1,\dots, u_N.$ 
\end{itemize}    
A system $(u_1,\dots, u_N)$ of unit vectors is called an {\it equiangular tight frame} ETF$(m,N)$ if it is tight and
the value of $|\langle u_i, u_j\rangle|$ is constant over all $i\neq j.$
It is well known (see e.g. \cite[Theorem 5.7]{FR} ) that if $u_1,\ldots u_N\in \mathbb{K}^m$ is an ETF then
\begin{equation}\label{Welch Bound}
\varphi_{m,N}:=|\langle u_i,u_j \rangle|=\sqrt{\frac{N-m}{m(N-1)}}
\qquad \textrm{ for all } i,j\in \{1,\dots, N \},\;i\neq j,
\end{equation}
and the cardinality $N$ cannot exceed $\frac{m(m+1)}{2}$ in the real case and $m^2$ in the complex case (see e.g. \cite[Theorem 5.10]{FR} ). An ETF that realizes this upper bound is called {\it the maximal ETF.}

\section{ Main Result.}
We start with a theorem, which is a special case of \cite[Lemma 5]{BC} that was stated for isotropic measures. In our case, the proof given by Bukh and Cox can be streamlined. 
\begin{thm}\label{BukhCox}
Let $1<m\leq N.$ Then the following inequalities holds
$$
\mu_\mathbb{R}(m,N)\leq \delta_{m,\frac{m(m+1)}{2}}=\frac{2}{m+1}\left(1+\frac{m-1}{2}\sqrt{m+2}\right)
$$
$$
\mu_\mathbb{C}(m,N)\leq \delta_{m,m^2}=\frac{1}{m}\left(1+(m-1)\sqrt{m+1}\right)
$$
\end{thm}
{\sc Proof.} Let $U\in\mathbb{K}^{m\times N}$ be such that $UU^*={\rm I}_m$. Denote by $u_i$  the i-th column of the matrix $U.$ Observe that for matrix $\widetilde{U}\in\mathbb{K}^{m\times \widetilde{N}}$ created from only nonzero columns (since $u_1,\dots, u_N$ form a tight frame, they span $\mathbb{K}^m$ and $\widetilde{N}\geq m$) we have
$$
\frac{1}{N}\sum_{i,j=1}^{N}|U^* U|_{ij}\leq \frac{1}{\widetilde{N}}\sum_{i,j=1}^{\widetilde{N}}|\widetilde{U}^* \widetilde{U}|_{ij}  \;\textrm{ and }\; \widetilde{U}\widetilde{U}^*={\rm I}_m\; .
$$
So, without loss of generality, we can assume that all $u_i$ are nonzero vectors.
Now, we are ready to define the matrices associated with the vectors $u_i$ by
\begin{equation}
    L_i=\|u_i\|^{-\frac{3}{2}}u_i u_i^*.
\end{equation}
Denote by $G$  the Gram matrix of the system $(L_1,\dots,L_N)$ with respect to the Frobenius inner product. Then,
\begin{align*}
G_{ij} &= \langle L_i, L_j \rangle_F = {\rm tr} (L_i L_j^*)= {\rm tr} (L_i L_j)= \|u_i\|^{-\frac{3}{2}}\|u_j\|^{-\frac{3}{2}}{\rm tr} (u_iu_i^{*}u_ju_j^*)\\
&= \|u_i\|^{-\frac{3}{2}}\|u_j\|^{-\frac{3}{2}}  |\langle u_i,u_j \rangle|^2
\end{align*}
 Now observe that for all $i\in\{1,\dots,N\}$, since $(u_1,\dots, u_N)$ is a tight frame with the constant $\alpha$ equal to 1
$$
\|u_i\|^{1/2}=\|u_i\|^{-\frac{3}{2}}\|u_i\|^2=\sum_{k=1}^N\|u_i\|^{-\frac{3}{2}} |\langle u_i, u_k \rangle|^2=\sum_{k=1}^N\|u_k\|^{\frac{3}{2}} G_{ik}.
$$
Let $g_k$ denote the  k-th column of the matrix $G$, then for fixed $i\in\{1,\dots,N\},$ by applying the latter we obtain
\begin{align*}
\sum_{k=1}^{N}\|u_i\|^\frac{1}{2}\|u_k\|^\frac{3}{2}g_k =[\|u_i\|^{\frac{1}{2}}\|u_1\|^{\frac{1}{2}},\dots,\|u_i\|^{\frac{1}{2}}\|u_N\|^{\frac{1}{2}}]^{\top}.
\end{align*}
Therefore, for any constants $a,b\in\mathbb{R}$ the matrix $A$ with the entries given by
$A_{ij}:=aG_{ij}-b\|u_i\|^{\frac{1}{2}}\|u_j\|^{\frac{1}{2}}$ has the rank less than or equal to ${\rm rk}(G).$ Let $\varphi>0.$ Then combining the Cauchy - Schwarz inequality and  ${\rm tr}(AA^*)\geq\frac{({\rm tr}(A))^2}{{\rm rk}(A)}$ (see e.g. \cite[Fact 7.12.13]{Be})  we get
\begin{align*}
\sum_{i,j=1}^N &\frac{(|\langle u_i,u_j\rangle|-\varphi\|u_i\|\|u_j\|)^2}{\|u_i\|\|u_j\|}=\sum_{i,j=1}^N \frac{(|\langle u_i,u_j\rangle|^2-\varphi^2\|u_i\|^2\|u_j\|^2)^2}{\|u_i\|\|u_j\|(|\langle u_i,u_j\rangle|+\varphi\|u_i\|\|u_j\|)^2}\\
&\geq \sum_{i,j=1}^N \frac{(|\langle u_i,u_j\rangle|^2-\varphi^2\|u_i\|^2\|u_j\|^2)^2}{(1+\varphi)^2\|u_i\|^3\|u_j\|^3}= \sum_{i,j=1}^N \left(\frac{G_{ij}}{1+\varphi}-\frac{\varphi^2}{1+\varphi}\|u_i\|^{\frac{1}{2}}\|u_j\|^{\frac{1}{2}}\right)^2\\
&\geq \frac{1}{{\rm rk}(G)}\left((1-\varphi)\sum_{i=1}^N\|u_i\|\right)^2 
\end{align*}
Squaring the addends of the left-sided sum and rearranging the latter inequality gives 
\begin{equation}\label{inq2}
    2\varphi\sum_{i,j=1}^{N}|U^*U|_{ij}\leq\sum_{i,j=1}^{N}\frac{|\langle u_i, u_j\rangle|^2}{\|u_i\|\|u_j\|} +\left(\varphi^2-\frac{(1-\varphi)^2}{{\rm rk}(G)}\right)\left(\sum_{i=1}^{N}\|u_i\|\right)^2.
\end{equation}
Now observe that using the Cauchy - Schwarz inequality and the tightness of the vectors $u_1,\dots,u_N$, we have
\begin{align*}
    \sum_{i,j=1}^{N}\frac{|\langle u_i, u_j\rangle|^2}{\|u_i\|\|u_j\|}&\leq \sqrt{\sum_{i,j=1}^{N}\frac{|\langle u_i, u_j\rangle|^2}{\|u_i\|^2}}\sqrt{\sum_{i,j=1}^{N}\frac{|\langle u_i, u_j\rangle|^2}{\|u_j\|^2}}=\sum_{i=1}^N\frac{1}{\|u_i\|^2}\sum_{j=1}^{N}|\langle u_i,u_j\rangle|^2\\
   & =\sum_{i=1}^N\frac{\|u_i\|^2}{\|u_i\|^2}=N
\end{align*}
and
\begin{align*}
\left(\sum_{i=1}^{N}\|u_i\|\right)^2=\langle[1,\dots,1], [\|u_1\|,\dots, \|u_N\|]\rangle^2 &\leq N \sum_{i=1}^{N}\|u_i\|^2 = N\;{\rm tr}(UU^*)\\ 
&=N\;{\rm tr}(U^*U)=N\;m.
\end{align*}
In view of the above, if we take $\varphi>0,$ such that $\varphi^2-\frac{(1-\varphi)^2}{{\rm rk}(G)}\geq 0,$ \eqref{inq2} reads
\begin{equation}\label{inq2u}
    2\varphi\sum_{i,j=1}^{N}|U^*U|_{ij}\leq N +\left(\varphi^2-\frac{(1-\varphi)^2}{{\rm rk}(G)}\right)N\; m.
\end{equation}
The rank of the Gram matrix of the system of vectors is equal to the dimension of the space spanned by these vectors. Now, we have to consider the real and complex cases separately. When $\mathbb{K}=\mathbb{R},$ then the matrices $L_1,\dots,L_N$ live in the space of symmetric matrices, whose dimension is equal to $\frac{m(m+1)}{2}.$ Therefore, ${\rm rk}(G)\leq \frac{m(m+1)}{2}.$ Setting $\varphi=\frac{1}{\sqrt{m+2}}$ (notice that if ETF$(m, \tfrac{m(m+1)}{2})$ exists $\varphi=\displaystyle\varphi_{m, \frac{m(m+1)}{2}}$),  \eqref{inq2u} turns into

\begin{equation}\label{ineq1}
\frac{1}{N}\sum_{i,i=1}^{N}|U^*U|_{ij}\leq\frac{2}{m+1}\left(1+\frac{m-1}{2}\sqrt{m+2}\right),
\end{equation}
which after taking the maximum over all $U$ leads to the announced upper bound of $\mu_{\mathbb{R}}(m, N).$\\ 
As the set of Hermitian matrices is not a vector space, in the complex case we only know that ${\rm rk}(G)\leq m^2.$ Now, fixing $\varphi=\frac{1}{\sqrt{m+1}}$ (if $ETF(m, m^2)$ exists, then $\varphi=\varphi_{m,m^2}$) and reasoning as in the real case, we obtain
$$
\mu_\mathbb{C}(m,N)\leq \delta_{m,m^2}=\frac{1}{m}\left(1+(m-1)\sqrt{m+1}\right)
$$
as desired.\\

\noindent Now to establish the upper bound for $\lambda_{\mathbb{K}}(m)$ it is sufficient to prove the equality of the maximal and quasimaximal absolute projection constants. In the real case, Basso proved this in \cite[Proof of Theorem 1.2]{B}. We propose an alternative argument, which is also valid in the complex case. 

\begin{thm}\label{basso}
Let $m\geq 1,$ then
\begin{equation}
    \lambda_\mathbb{K}(m)=\mu_\mathbb{K}(m).
\end{equation}
\end{thm}
{\sc Proof.} Since for every $N>m$ we have $\mu_\mathbb{K}(m,N)\leq\lambda_\mathbb{K}(m,N)$, it is enough to show that $\lambda_\mathbb{K}(m,N)\leq \mu_{\mathbb{K}}(m).$ We choose $U_0\in\mathbb{K}^{m\times N}$ and $t^0\in\mathbb{R}_{+}^n$ that realize $\lambda_\mathbb{K}(m,N).$ Next we define the function
$$
f: \mathbb{R}^N_+\ni t \mapsto \sum_{i,j=1}^{N}t_it_j|U_0^*U_0|_{ij}.
$$
Since $f$ is a continuous function, for every $\varepsilon>0$ there exists $t^\varepsilon \in \mathbb{Q}_+^N$ such that $\lambda_\mathbb{K}(m,N)=f(t^0)\leq f(t^\varepsilon)+\varepsilon$ and $\|t_0-t_\varepsilon\|\leq \varepsilon.$ Taking the common denominator, we know that there exists $q,n_1,\dots,n_N\in \mathbb{N}$ such that $t^\varepsilon=\frac{1}{q}(n_1,\dots,n_N).$ Observe that $\|t^\varepsilon\|q=\sqrt{n_1^2+\dots+n_N^2.}$ Put $\widetilde{N}=n_1^2+\dots+n_N^2.$ Let $u_i$ denote the i-th column of the matrix $U_0$. Now we consider the matrix $U_\varepsilon \in\mathbb{R}^{m\times \widetilde{N}} $ defined, in block notation, by 
$$
U_\varepsilon:=\left[\tfrac{1}{n_1}u_1\mathbb{1}_{n_1^2}^*|\dots|\tfrac{1}{n_N}u_N\mathbb{1}_{n_N^2}^*\right],
$$
where $\mathbb{1}_{k}$ denotes the k-dimensional vector with all entries equal to $1$. 
Notice that 
\begin{align*}
    U_\varepsilon U_\varepsilon^*&=\frac{1}{n_1^2} u_1\mathbb{1}_{n_1^2}^*\mathbb{1}_{n_1^2}{u_1}^*+\dots + \frac{1}{n_N^2} u_N\mathbb{1}_{n_N^2}^*\mathbb{1}_{n_N^2}{u_N}^*=u_1u_1^*+\dots+u_Nu_N^*\\
    &=U_0 U_0^*={\rm I_m}.
\end{align*}
In turn,
\begin{align*}
 \frac{1}{\widetilde{N}}\sum_{i,j=1}^{\widetilde{N}}|U_\varepsilon^* U_\varepsilon|_{i j}=  \frac{1}{\widetilde{N}}\sum_{i,j=1}^{N}\left|\left\langle \frac{1}{n_i}u_i,\frac{1}{n_j}u_j\right\rangle\right|n_i^2 n_j^2= \frac{1}
 {\| t^{\varepsilon}\|^2}\sum_{i,j=1}^{N}\left|\left\langle u_i,u_j\right\rangle\right|\frac{n_i}{q}\frac{n_j}{q}= \frac{f(t_\varepsilon)}{\|t^{\varepsilon}\|^2}.
\end{align*}
Putting everything together yields
$$
\mu_\mathbb{K}(m)\geq \mu_\mathbb{K}(m, \widetilde{N})\geq \frac{ \lambda_\mathbb{K}(m, N) - \varepsilon}{(1+\varepsilon)^2}.
$$
Letting $\varepsilon$ to $0$ we establish the desired result.
~\\
~\\
In view of Theorem \ref{SF} and Theorem \ref{BukhCox}, the latter leads immediately to the following corollary.  
\begin{thm}\label{MPC}
Let $ m>1$ then
\begin{enumerate}[i)]
    \item $\lambda_\mathbb{R}(m) \leq \delta_{m,\frac{m(m+1)}{2}}=\frac{2}{m+1}\left(1+\frac{m-1}{2}\sqrt{m+2}\right)\,$ 
    \item $ \lambda_\mathbb{C}(m)\leq \delta_{m,m^2}=\frac{1}{m}\left(1+(m-1)\sqrt{m+1}\right).$ 
\end{enumerate}
Furthermore, the inequality becomes an equality if there exists a maximal $ETF$ in $\mathbb{K}^m.$ 
\end{thm}
The above was first stated in \cite{KT}, with the proof based on an erroneous lemma,  as pointed out in \cite{CLe}.
In the real case, maximal ETFs seem to be rare objects. The only known cases are for $m$ equal to $2,$ $3,$ $7$ and $23.$ Applying the {\it Theorem \ref{MPC}}  we get
\begin{thm}
~
\begin{itemize}
    \item $\lambda_\mathbb{R}(2)=\frac{4}{3}.$ 
    \item $\lambda_\mathbb{R}(3)=\frac{1+\sqrt{5}}{2};$ 
    \item $\lambda_\mathbb{R}(7)=\frac{5}{2};$
    \item $\lambda_\mathbb{R}(23)=\frac{14}{3}.$
    \end{itemize}
    \end{thm}
It is worth mentioning that the first bullet point is known as the Gr\"unbaum conjecture, which was stated in \cite{G} and the first complete proof of it was given in \cite{CL}. Recently, a new proof was presented in \cite{B}. However, both proofs are long and quite complicated.\\ 
Unlike the real case, numerous examples of complex maximal ETF are known. For example, for $m \in\{1,\ldots, 17,19,24,28,35,48\}$ (see, e.g. \cite{FM}). Hence, we have the following.
\begin{thm}
For $m \in\{1,\ldots, 17,19,24,28,35,48\}$  we have
$$
\lambda_\mathbb{C}(m)=\frac{1}{m}\left(1+(m-1)\sqrt{m+1}\right).
$$
\end{thm}
In fact, it is conjectured that there is a complex maximal ETF in every dimension (Zaurner's conjecture \cite{Z}), which allows us to state the following conjecture.  

\begin{conjecture}
For every $m\geq1$
$$
\lambda_\mathbb{C}(m)=\frac{1}{m}\left(1+(m-1)\sqrt{m+1}\right).
$$
\end{conjecture}

\end{document}